\documentclass[12pt]{amsart}
\usepackage{epsfig}

\newcommand{\G}{Sl_2{\mathbb{C}}}

\newcommand{\rep}{\mathrm{Rep}}
\newcommand{\tor}{\pi_1(T^2)}

\newcommand{\ql}{\mathbb{C}_t[l,l^{-1},m,m^{-1}]}
\newcommand{\qp}{\mathbb{C}_t[l,m]}
\newcommand{\lcr}{\raisebox{-5pt}{\mbox{}\hspace{1pt}
                  \epsfig{file=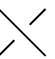}\hspace{1pt}\mbox{}}}
\newcommand{\ift}{\raisebox{-5pt}{\mbox{}\hspace{1pt}
                  \epsfig{file=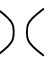}\hspace{1pt}\mbox{}}}
\newcommand{\zer}{\raisebox{-5pt}{\mbox{}\hspace{1pt}
                  \epsfig{file=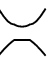}\hspace{1pt}\mbox{}}}

\title{The A-Polynomial  from the Noncommutative Viewpoint}
\author{Charles Frohman}
\address{Department of Mathematics, University of Iowa, Iowa City, IA
52242, USA}
\email{frohman@math.uiowa.edu}

\author{R{\u{a}}zvan Gelca}
\address{Department of Mathematics,University of Michigan, Ann Arbor, MI 48109, USA}
\email{rgelca@math.lsa.umich.edu}
\author{Walter Lofaro}
\address{Department of Mathematics, Boise State University,
 Boise, Id 83725, USA}
\email{lofaro@math.idbsu.edu}

\newtheorem{theorem}{Theorem}
\newtheorem{prop}{Proposition}

\begin{document}
\maketitle
\section{Introduction}
This paper places the A-polynomial of a knot into the framework of noncommutative geometry.
The A-polynomial was introduced in \cite{CCGLS}. The A-polynomial describes how the
$\G$-characters of a knot lie inside the $\G$-characters of its boundary torus.
It is known that the A-polynomial can be related to the Alexander polynomial of
the knot and to the structure of essential  surfaces in the complement of the knot
\cite{CL}.

We construct an invariant depending on a complex parameter, by following
the construction of the A-polynomial, but replacing 
$\G$-characters of surfaces and 3-manifolds by the
Kauffman bracket skein module of cylinders over surfaces and 3-manifolds. The
Kauffman bracket skein module depends on a parameter $t$, so that when $t$ is
set equal to $-1$ the $\G$-characters are recovered. The A-polynomial can be derived from the 
noncommutative invariant when $t$ is set equal to $-1$. 

The construction of the noncommutative invariant uses the relationship between
the Kauffman bracket skein module and the noncommutative torus, already explicated
in \cite{FG}. The noncommutative torus is one of the fundamental examples of a noncommutative
space \cite{Co}. 

The invariant is a minimal reduced Gr\"{o}bner basis of a left ideal in the noncommutative
plane. The use of Gr\"{o}bner bases is an important theme in modern commutative algebra.
For instance, algorithms associated with Gr\"{o}bner bases are the heart of most
symbolic manipulation programs. The use of Gr\"{o}bner bases to study the noncommutative
plane is foreshadowed by the work of Cohn on skew polynomial rings \cite{C}.

There is a left action of the Kauffman bracket skein module of a cylinder over a torus
on the Kauffman bracket skein module of a solid torus. This allows us to compare the
noncommutative invariant with data coming from colored Jones polynomials. Specifically, the
noncommutative invariant annihilates the data from the colored Jones polynomials. This should 
lead to a deeper understanding of the relationship between
the Jones polynomial and the representation theory of the fundamental group of the complement
of a knot.

In section 2 we recall the basic definitions associated with the A-polynomial. In section 3
we introduce the noncommutative analogues of the spaces and maps used in section 2 and then
define the noncommutative A-ideal and the A-basis. In section 4 we derive the orthogonality relationship
between the noncommutative A-ideal and the Jones polynomial. In the final section we pose some
questions and suggest avenues of future investigation.

\section{Characters and Knot Invariants}

 \subsection{ $\G$-Characters } 

Let $G=<a_i | r_j>$ be a finitely generated group. 
A representation $\rho:G \rightarrow \G$ is determined by a choice of matrices $A_i$ in $\G$ 
so that
when the relations $r_j$ are rewritten with the $a_i$ replaced with the $A_i$
they are equal to the identity in $\G$. In other words, the representations
can be identified with the subset of $\prod_{i=1}^n \G$, where $n$ is the
number of generators of $G$, that satisfies the equations obtained by requiring the $r_j$ to
evaluate to the identity.
Denote this subset $\mathrm{Rep}(G)$, and call it the representations of $G$.
Since we only deal with $\G$-representations, we suppress references to $\G$.

Let $a_{rs}(i)$, $r,s \in \{1,2\}$, and $i \in \{1,\ldots ,n\}$ be the function on
 $\prod_{i=1}^n \G$ which yields
the entry in the r-th row and
s-th column of the i-th entry of each element of  $\prod_{i=1}^n \G$. The
coordinate ring $C[\prod_{i=1}^n \G]$ consists of the 
 polynomials in the $a_{rs}(i)$ modulo the ideal
generated by 
\[a_{11}(i)a_{22}(i)-a_{12}(i)a_{21}(i)-1,\] 
for $i \in \{1,\ldots ,n\}$.
Each relation $r_j$ gives rise to four polynomials corresponding to the four entries of a matrix,
and arising from the condition that the relations evaluate to  the identity.
Let $I(G)$ be the ideal in $C[\prod_{i=1}^n \G]$
generated by the polynomials coming from the relations. 
Finally, $R(G)=C[\prod_{i=1}^n \G]/I(G)$
is the {\em affine representation ring} of $G$. In more generality, Lubotsky and Magid \cite{L-M}
proved that the isomorphism class of $R(G)$ is an invariant of the group. The word
``affine'' refers to the fact that $R(G)$ is the unreduced coordinate ring of the
representations. To get the classical representation ring, take the quotient
of $R(G)$ by its nilradical $\sqrt{0}=\{p | p^n=0 \  \mathrm{ for \  some} \  n \}$.

There is a left action of $\G$ on $\prod_{i=1}^n \G$ by conjugation,
\[ A\bullet (A_1, \ldots , A_n)=(AA_1A^{-1}, \ldots ,AA_nA^{-1}).\]
This action induces a right action on $C[\prod_{i=1}^n \G]$. It is easy
to check that the action leaves $I(G)$ invariant, hence the action descends
to a right action on $R(G)$. The invariant subring of this action,
call it $\chi(G)=R(G)^{\G}$,
is called the {\em affine characters} of $G$. This ring is an invariant of the group $G$.
Once again, to obtain what is classically referred to as the characters, take
the quotient of $\chi(G)$ by its nilradical.

For the most frequently encountered groups we are being too careful, the affine
representation ring and the affine character ring have trivial nilradical. However,
there are examples of groups where the distinction is real \cite{KM}. 

The affine characters have recently been the subject of scrutiny in works of
Bullock, Brumfiel-Hilden, Przytycki-Sikora, and Sikora, \cite{B,BH,PS,S}. 
Here is an intrinsic
definition due to Sikora. Let $\hat{G}$ be the set of conjugacy classes of $G$. If $W$ is
an element of $G$ we denote the conjugacy class of $W$ by $<W>$.  Let $S(\hat{G})$ be the
symmetric algebra on  $\hat{G}$, that is,  polynomials where the variables are conjugacy
classes in $G$ and the coefficients are complex numbers.   Let $J$ be the ideal generated 
by all polynomials , $<Id>+2$,
and $<AB>+<A^{-1}B>+<A><B>$ where $A$ and $B$ range over the elements of $G$. There
is an isomorphism $S(G)/J \rightarrow \chi(G)$ induced by sending each $<A>$ into
the polynomial corresponding to $-\mathrm{tr}(A)$.

There is a space corresponding to the characters. The action of
$\G$ on $C[\prod_{i=1}^n \G]$ is not good, in the sense that the quotient is not
Hausdorff. We sidestep this by defining an equivalence relation on $\rep (G)$, that
yields a Hausdorff space. Let  $\rho, \eta \in \rep (G)$ be equivalent
if for every $t \in \chi(G)$, $t(\rho)=t(\eta)$. The resulting quotient space
is the character variety of $G$. An {\em algebraic subset} of an 
$\mathbb{C}^n$  is the 
solution set of a system of polynomial equations. A set of polynomials $S$ {\em cuts out} an 
algebraic subset
$V$, if the intersection of the zeroes of elements of $S$ is exactly $V$. The ideal $I(V)$
of $V$ is the set of all polynomials that vanish on $V$. It is a consequence
of the Nullstellensatz that the ideal of $V$ is the radical of the smallest ideal containing
any set of polynomials that cut $V$ out. The coordinate ring of $V$
is the quotient of the ring of polynomials in $n$ variables by $I(V)$. Sometimes
it is nice to have a set $X(G)$ corresponding to the characters $\chi(G)$. Of course 
the set only
matters to the extent that its points correspond to maximal ideals in $\chi(G)$. To
this end, define 
a {\em realization} of $X(G)$ to be an algebraic
subset of some $\mathbb{C}^n$, where the coordinates on $\mathbb{C}^n$ correspond to traces of elements 
of $G$, whose coordinate  ring is isomorphic to $\chi(G)/\sqrt{0}$.
The set $X(G)$ is in one to one correspondence with the character variety of $G$  \cite{CS}.

\subsection{Characters of the Torus} 

Concentrate on the case of the 
fundamental group of a torus.
Let us think of it as the free Abelian group on $\lambda$ and $\mu$. The representations
correspond exactly to pairs of matrices in $\G$ that commute. Next we describe three
functions on the space of representations.  Let $x$ be the trace of the image of $\lambda$,
$y$ be the trace of the image of $\mu$ and let $z$ be the trace of the image of their product.
The affine character ring, $\chi(\pi_1(T^2))$, is
generated by $x$, $y$, $z$ with relation 
\[x^2+y^2+z^2-xyz -4=0. \]
The ring has
trivial nilradical. The space $X(\tor)$ is realized as the set of points in $\mathbb{C}^3$
satisfying the equation above.  There is a two-fold branched cover
of $X(\tor)$ that is used in the definition of knot invariants. Consider 
$\mathbb{C}^* \times \mathbb{C}^*$, where $\mathbb{C}^*$ denotes the nonzero complex numbers, give it coordinates $l$ and $m$. In order to see that  its coordinate ring 
is $\mathbb{C}[l,l^{-1},m,m^{-1}]$, it is helpful
to think of $\mathbb{C}^* \times \mathbb{C}^*$ as ordered pairs of diagonal matrices
in $\G$.
\[ \left(\begin{pmatrix} l & 0 \\ 0 & l^{-1} \end{pmatrix},
\begin{pmatrix} m & 0 \\ 0 & m^{-1} \end{pmatrix}\right)\]

There is a map $p: \mathbb{C}^* \times \mathbb{C}^*  \rightarrow X(\tor)$ given
by sending each pair of points to its equivalence class. The map $p$ is a two-fold 
branched cover whose singular points are the four ordered pairs chosen from
$\begin{pmatrix} \pm 1 & 0 \\ 0 & \pm 1 \end{pmatrix}$.  In terms of our realization
of $X(\tor)$, the map sends 
\[ \left(\begin{pmatrix} l & 0 \\ 0 & l^{-1} \end{pmatrix},
\begin{pmatrix} m & 0 \\ 0 & m^{-1} \end{pmatrix}\right)\]
to the triple $(l+l^{-1}, m+ m^{-1}, lm +l^{-1}m^{-1})$. Being a two-fold
branched cover, $p$ has a deck transformation, 
\[\theta:  \mathbb{C}^* \times \mathbb{C}^* \rightarrow  \mathbb{C}^* \times \mathbb{C}^*,\]
$\theta(l,m)=(l^{-1},m^{-1})$.

Dual to $p$, and $\theta$ are maps, that we  denote $\hat{p}$ and $\hat{\theta}$, 
\[\hat{p}: C[X(\tor)] \rightarrow  \mathbb{C}[l,l^{-1},m,m^{-1}],\]
and 
\[ \hat{\theta}:  \mathbb{C}[l,l^{-1},m,m^{-1}] \rightarrow \mathbb{C}[l,l^{-1},m,m^{-1}].\]
One can show that the image of $\hat{p}$ is the fixed subalgebra of $\hat{\theta}$.

\begin{prop} Suppose that $V \subset X(\tor)$ is algebraic and $S$ is a set
 of functions that cuts out $V$. The ideal of  $p^{-1}(V)$
is the radical of the smallest ideal of $ \mathbb{C}[l,l^{-1},m,m^{-1}]$ containing
$\hat{p}(S)$.\end{prop}

\proof The set $\hat{p}(S)$ cuts out $p^{-1}(V)$, the rest is a standard characterization
of the ideal of an algebraic set. \qed

The coordinate ring $\mathbb{C}[l,l^{-1},m,m^{-1}]$ can be seen as the ring of
fractions of $\mathbb{C}[l,m]$ with respect to the set of monomials in $l$ and $m$.

\begin{prop}If $I \subset \mathbb{C}[l,l^{-1},m,m^{-1}]$ is an ideal then its contraction
$J$ to $\mathbb{C}[l,m]$ is $I \cap \mathbb{C}[l,m] $. The extension of $J$
is just the ideal of  $\mathbb{C}[l,l^{-1},m,m^{-1}]$ generated by $J$.
Furthermore, the extension of $J$ is $I$. \end{prop}
 \qed

\subsection{The A-polynomial} The A-polynomial is defined in \cite{CCGLS}. Let
$K$ be a knot in $S^3$. Let $\mathrm{Rep}_D$ denote the set of representations of the
fundamental group of the complement of $K$ in which the peripheral subgroup is sent to
diagonal matrices. Let $\rho : Rep_D \rightarrow \mathbb{C}^* \times \mathbb{C}^*$
be the map that sends each representation to the ordered pair consisting of the upper
diagonal entry of the images of the longitude and meridian of the knot. The image
of $\rho$ is an algebraic curve. Therefore, its ideal is principal and generated
by a Laurent polynomial $B(l,m)$.  It can be shown that the image of $\rho$ always has
the set $l=1$ as a component, so $B(l,m)=(l-1)A(l,m)$, where $l-1$ does not divide $A(l,m)$.
Properly normalized, the
polynomial $A(l,m)$ is an invariant of $K$, called {\em the A-polynomial}.
It is clear that $A(l,m)$ can be recovered from $B(l,m)$.

There is a map $r:\mathrm{Rep}(\pi_1(S^3-K)))  \rightarrow \mathrm{Rep}(\tor)$
given by restriction. This induces a map 
\[ \hat{r}: \chi(\tor)) \rightarrow \chi(\pi_1(S^3-K)) .\]
The kernel of $\hat{r}$ is an ideal $\mathcal{I}(K)$. Notice $\mathcal{I}(K)$  cuts out an algebraic 
subset $V$ of $X(\tor)$.

\begin{prop} The radical of the extension of   $\hat{p}(\mathcal{I}(K))$ in 
$\mathbb{C}[l,l^{-1},m,m^{-1}]$ is the principal ideal generated by $B(l,m)$.
That is the ideal $\mathcal{I}(K)$ determines the A-polynomial. \end{prop}

\proof The ideal $\mathcal{I}(K)$ is exactly the functions that vanish on the image of $r$.  
>From the definitions,  $p^{-1}(\mathrm{im}(r))$ is equal to the image of $\rho$. The first
proposition then finishes the proof. \qed

You can think of the extension of the  ideal $\mathcal{I}(K)$ as the holomorphic sections of a line
bundle over $\mathbb{C}^* \times \mathbb{C}^*$. Specifically, the sections of
the line bundle associated
to the divisor of $\frac{1}{B(l,m)}$.

\section{Noncommutative Characters and Knot Invariants}

\subsection{The Noncommutative Torus} 




The noncommutative torus is 
a ``virtual'' geometric space whose algebra of 
continuous functions is the (irrational) rotation
algebra $A_\theta$ \cite{Ri}. It is customary to call the
algebra of functions itself the noncommutative torus.

The algebra $A_\theta$ is usually defined for a real angle of
rotation $\theta$, however we consider 
$\theta$ to be any complex number, and let $t=e^{\pi i\theta}$.
This 
algebra can be introduced abstractly by exponentiating
the Heisenberg noncommutation relation. That is,
$A_{\theta}$ is the closure in a certain $C^*$-norm of the
algebra spanned by $l,m,l^{-1},m^{-1}$, subject to the relation 
$lm=t^2ml$.

Rieffel \cite{Ri} defined $A_{\theta}$ in a concrete setting. Let 
$e_{p,q}= t^{-pq} l^pm^q$, where $p,q \in \mathbb{Z}$.
He defines his multiplication via the formula

\[
e_{p,q}*e_{r,s}=t^{|^{pq}_{rs}|}e_{p+r,q+s}
\]
which, from our approach, is just a consequence of the defining relation.
He  then considered the
closure of this algebra in the norm determined by the left
regular representation on $L^2(T^2)$. 

For the purpose of this paper we are interested only in the subalgebra
of the noncommutative torus consisting of
Laurent polynomials in $l$ and $m$, which we denote by
$\mathbb{C}_t[l,l^{-1},m,m^{-1}]$.
There is an  automorphism 
\begin{eqnarray*}
{\Theta} :\ql\rightarrow\ql, \; {\Theta} (e_{p,q})=e_{-p,-q}
\end{eqnarray*}
Let $\ql^{{\Theta}}$ be its invariant part
(which is spanned by $e_{p,q}+e_{-p,-q}$, $p,q\in{\mathbb{Z}}$, i.e.
by the noncommutative cosines).

In addition, let $ \mathbb{C}_t[l,m]$ be the subalgebra of $\ql$ spanned by
$e_{p,q}$, with $p,q \geq 0$. This is nothing but the ring of  noncommuting 
polynomials
in two variables $l$ and $m$  satisfying the noncommutation 
relation $lm=t^2ml$. This ring is frequently referred to as the noncommutative
plane.

\begin{prop}(\cite{Kas}, Proposition
IV.1.1) The ring $\qp$ is left Noetherian and  has no zero divisors. 
\end{prop}


We now need to broach the subject of Gr\"{o}bner bases in $\qp$.
As $\qp$ is so close to being commutative, the concepts translate
over very easily from the case of two variable polynomials.

We lexicographically order the $l^pm^q$. Hence $l^pm^q < l^rm^s$ if
either $p < r$ or $p=r$ and $q < s$. Given $f \in \qp$ we can write
$f= \sum \alpha_{p,q} e_{p,q}$ where the sum is finite. The leading term
$lt(f)$
of $f$ is the $\alpha_{p,q} l^pm^q$ where the $l^pm^q$ is largest in
the lexicographical ordering among those terms with $\alpha_{p,q} \neq 0$.
The leading power product is $l^pm^q$ and the leading coefficient is
$\alpha_{p,q}$.

Suppose that $u,v,w \in \qp$ and $u=vw$, then we say $w$ divides $u$ on the right
and we let $\frac{u}{w}=v$. A {\em Gr\"{o}bner basis} for a left ideal $I$ is a collection
$f_i$ of elements of the ideal $I$ so that the ideal generated by the leading terms
of the $f_i$ is equal to the ideal generated by the leading terms of elements of $I$.
We say the Gr\"{o}bner basis is {\em minimal} if no two $f_i$ have the same leading power product.
We say the Gr\"{o}bner basis is {\em reduced} if no power product in each $f_i$ is divisible
by the leading power product of any other $f_i$. 

\begin{prop} 
For an ideal of polynomials in $\qp$ there exists a unique minimal, reduced 
 Gr\"{o}bner basis,  consisting  of monic polynomials.
 \end{prop}

\proof  By changing any
statements in \cite{AL} about ideals to statements about left ideals, the proof
there goes through verbatim. \qed

\subsection{ The Kauffman Bracket Skein Module}

Let $M$ be an orientable $3$-manifold.
A framed link in $M$ is an embedding of a disjoint union of annuli
into $M$. In diagrams we will draw only the core of an annulus lying
parallel to the plane of the paper (i.e. with blackboard framing).

Two framed links in $M$ are equivalent if there is an isotopy of $M$
taking one to the other. Let $\mathcal{L}$ denote the set of equivalence
classes of framed links in $M$, including the empty link. Fix a complex number
$t$.  Consider the vector space,  $\mathbb{C} \mathcal{L}$ with  basis $\mathcal{L}$.
Define $S(M)$ to be the smallest subspace of $\mathbb{C} \mathcal{L}$
containing all expressions of the form
$\displaystyle{\lcr-t\zer-t^{-1}\ift}$
and 
$\bigcirc+t^2+t^{-2}$,
where the framed links in each expression are identical outside 
balls, in which they look like
 pictured in the diagrams. The Kauffman bracket skein module $K_t(M)$ is
the quotient
\[ \mathbb{C} \mathcal{L} / S(M). \]
In the case of the cylinder over the torus, $T^2\times I$,  $K_t(T^2 \times I)$ 
has the structure of an algebra with multiplication given by
laying one link over the other. More precisely, to multiply skeins
corresponding to links $\alpha$ and $\beta$, isotope them so that
$\alpha$ lies in $T^2\times [\frac{1}{2},1]$ and $\beta$ in 
$T^2\times [0,\frac{1}{2}]$. Then $\alpha\cdot\beta$ is the element of
the skein module represented by the class of the union of these two
links in $T^2\times [0,1]$. Extend this to a distributive product.

Oriented simple closed curves on the torus up to isotopy are indexed by pairs of relatively prime
integers $(p,q)$. Corresponding to $(p,q)$  is a framed link in $T^2 \times I$.
Take  an annulus in $T^2 \times I$ whose core projects to a  $(p,q)$ curve, so that the
annulus runs parallel to the boundary of $T^2 \times I$. As the framed links are unoriented,
$(p,q)$ and $(-p,-q)$ give rise to the same link. Denote
the link corresponding to the $(p,q)$ by $L_{p,q}$. A standard argument based
on the proof that the Kauffman bracket in $S^3$ is well defined shows that as a vector
space, $K_t(T^2 \times I)$ has as basis all links consisting of parallel copies of the
$L_{p,q}$.

Let $x$ be  $L_{0,1}$ , $y=L_{1,0}$ and $z=L_{1,1}$. It is a theorem of Bullock and
Przytycki \cite{BP} that $K_t(T^2 \times I)$ is isomorphic to polynomials in three noncommutative
variables, $x$, $y$ and $z$ modulo the ideal generated by

\[ t^2x^2 + t^{-2}y^2 + t^2z^2 -txyz-2(t^2 + t^{-2}), \]
\[ txy-t^{-1}yx -(t^2 -t^{-2})z,\]
\[ tzx-t^{-1}xz -(t^2 -t^{-2})y,\]
and
\[ tyz-t^{-1}zy -(t^2 -t^{-2})x.\]

When $t=-1$ the Kauffman bracket skein module of an arbitrary three-manifold
 can be made into an algebra.
The point is that at $t=-1$ the
skein relation allows us to change crossings. To multiply two links, perturb them so that they miss one another
and take their union. As crossings don't count, the answer is independent of the
perturbation. This extends to make $K_{-1}(M)$ into an algebra for any $M$. It is
a theorem of Bullock \cite{B} that $K_{-1}(M)/\sqrt{0}$ is naturally  isomorphic to
$\chi(\pi_1(M))/\sqrt{0}$ . You can see this isomorphism from the description of
$\chi(\pi_1(M))$ due to Sikora given above. The correspondence at $t=-1$ is slightly
tricky as the $x$, $y$ and $z$ given here correspond to $-x$, $-y$ and $-z$ in
the relation we gave for the $SL_2{\mathbb{C}}$-characters of the torus.

In \cite{FG} the following theorem is proved.

\begin{theorem} There exists an isomorphism of
algebras 
\[\hat{p} : K_t({ T}^2\times I)\rightarrow \ql^{\Theta}\]
determined by  

\[\hat{p}(L_{p,q})=e_{(p,q)}+e_{(-p,-q)}, \: p,q\in {\mathbb{Z}}.\]

\end{theorem}

\subsection{ The noncommutative A-basis}
Let $M$ be an oriented homology 3-sphere. Let $K$ be a knot in $M$. Let $X$ be the complement
of a regular neighborhood of $K$. The space $X$ is a compact manifold with boundary
homeomorphic to $T^2$. The meridian of the knot and longitude are well defined up
to sign, by requiring that the meridian bound a disk in the regular neighborhood
and the longitude bound a Seifert surface, and their intersection number in the boundary
of $X$ is 1. There is a map $\hat{r} : K_t(T^2 \times I) \rightarrow K_t(X)$ obtained
by gluing the cylinder over the torus into the complement of the knot at the $T^2 \times
\{ 0 \}$ end so that the meridian goes to the meridian and the longitude goes to the
longitude.  Let $\mathcal{I}_t(K)$ be the left ideal which is 
the kernel of $\hat{r}$. This ideal is called the peripheral ideal
of  the knot $K$. The map $\hat{p}$
takes $\mathcal{I}_t(K)$ to $\mathbb{C}_t[l,l^{-1},m,m^{-1}]$. Let $\mathcal{J}_t(K)$ be
the extension of $\mathcal{I}_t(K)$ to a left  ideal in $\ql$.
Finally, the noncommutative A-ideal  $\mathcal{A}(K)$ is  the contraction  of $\mathcal{J}_t(K)$ to
$\qp$.

The {\em A-basis} of a knot $K$ at the complex number $t$ is the reduced
minimal Gr\"{o}bner basis of the left ideal $\mathcal{A}_t(K)$.

Recall for any 3-manifold  $K_{-1}(M)/\sqrt{0}$ is isomorphic to $\chi(\pi_1(M))/\sqrt{0}$.
For a cylinder over a torus the radical of both rings is trivial, hence
\[  K_{-1}(T^2 \times I)=\chi(\tor).\]
It can be shown that under this isomorphism, the ideal $\mathcal{I}_{-1}(K)$ cuts out
the image of the characters of $\pi_1(S^3-K)$, and hence its extension cuts out
the variety used to determine the A-polynomial.

\begin{theorem} The noncommutative A-basis at $t=-1$ determines the A-polynomial \end{theorem}

\qed

Let us denote by ${\mathcal I}_t$ the left ideal that 
is the kernel of the epimorphism
\begin{eqnarray*}
\hat{\rho} :K_t({ T}^2\times [0,1]) \rightarrow K_t({ D}^2\times 
{ T}).
\end{eqnarray*}

In \cite{FG} it is proved that $\mathcal {I}_t$ is  generated by  $L_{(0,1)}+t^2+t^{-2}$ and
$L_{(1,1)}+t^{-3}L_{(1,0)}$. To compute the noncommutative A-basis one first
shows that the contraction of $\mathcal{I}_t$ is generated by
$m^2+(t^2+t^{-2})m+1$ and $t^{-3}l^2m^2+t^{-5}l^2m+t^{-1}m+t$.  These are just
the contractions of the generators of $\mathcal {I}_t$ to the noncommutative
plane.  Then an easy application of Buchberger's Algorithm yields the
minimal reduced Gr\"{o}bner basis $\{m^2+(t^2+t^{-2})m+1, l^2m+t^{-2}l^2-m-t^2\}$.
To get this to correspond to the value of the A-polynomial there are two subtleties.
The first is that as we set it up, the meridian bounds instead of the longitude.
The second is the negative sign arising in the correspondence between skeins 
and characters. If you exchange $l$ with $-m$ and $m$ with $-l$, then the least
common multiple of the two polynomials will be $l-1$ which is as one would expect.


\section{Relation with the Jones polynomial}
 Assume in this section that $t$ is not a root of unity. The definition and 
normalization of the Jones-Wenzl idempotents will be used as in \cite{Li}. Let $S_c$ be the skein
in the solid torus obtained by plugging the $c$-th Jones-Wenzl idempotent into the
core of the solid torus. The Kauffman bracket skein module of the solid torus
has the set  $\{S_c\}$ as a basis.  Let $\hat{K}_t(S^1\times D^2)$ be the vector
space of formal sums $\sum_c z_c S_c$ where the $z_c$ are complex numbers
and the $c$ range over the natural numbers starting at $0$.

The double of the solid torus is $S^1\times S^2$. Any skein $\alpha$ in $S^1\times S^2$
can be represented by a linear combination of framed links that miss $1 \times S^2$.
Hence, any skein in $S^1\times S^2$ can be represented  as a skein in a punctured ball, 
and  it has a Kauffman bracket.
 It follows that $K_t(S^1\times S^2)$ is canonically isomorphic to $\mathbb{C}$. There is 
a pairing
between skeins in $K_t(S^1 \times D^2)$. If $\alpha, \beta \in K_t(S^1 \times D^2)$
are represented by a single framed link each, take the union of two copies of the
solid torus, identified along their boundaries, with the link representing $\alpha$
in one and the link representing $\beta$ in the other. As this yields  a skein
in $S^1\times S^2$ we get a complex number by taking the Kauffman bracket as above. This can be extended
bilinearly to give a pairing,
\[  K_t(S^1 \times D^2) \otimes \hat{K}_t(S^1\times D^2) \rightarrow \mathbb{C},\]
for although the sum is infinite only finitely many terms are nonzero. In this
way we identify $ \hat{K}_t(S^1\times D^2)$ with $K_t(S^1\times D^2)^*$.

There is a representation of the Kauffman bracket skein algebra of the cylinder
over the torus into endomorphisms of $\hat{K}_t(S^1\times D^2)$. Glue the cylinder
onto $S^1 \times D^2$ along the $0$-end of the cylinder so that longitudes and
meridians go to longitudes and meridians. The 
matrix of a skein in the cylinder over the torus as a matrix with respect to the basis
$S_c$, is of bounded 
width. That is, there is an integer $n$ so that if $|i-j|>n$ the $ij$-entry
of the matrix is zero. 
 This can be seen,
as the matrix induced is the adjoint of the matrix corresponding to the
endomorphism of $K_t(S^1\times D^2)$ induced by gluing the $1$-end of the
cylinder to $S^1\times D^2$. 
 Notice that if $\hat{Z} \in \hat{K}_t(S^1\times D^2)$
then the annihilator of $\hat{Z}$ in the Kauffman bracket skein module of the
cylinder over the torus is a left  ideal.

Let $K \subset S^3$ be a framed knot. Let $X$ be the complement of an open
regular neighborhood of the knot.  There is a pairing,
\[ K_t(S^1\times D^2) \otimes K_t(X) \rightarrow \mathbb{C},\]
obtained by gluing the solid torus into the knot so that the meridian
of the solid torus goes to the meridian of the knot and the blackboard
longitude goes to the framing of the knot. To pair two skeins, take
their union and then take the Kauffman bracket in $S^3$ of the result. By using
the empty skein in $X$ we get a linear functional,
\[ Z(K) :  K_t(S^1\times D^2) \rightarrow \mathbb{C}.\]
 Let $\kappa(K,c)$ be the value of $Z(K)$ on $S_c$. We can then represent
$Z(K)$ by 
\[ \hat{Z}(K) = \sum_c \kappa(K,c)\sigma_c\]
where $\{\sigma _c\}$ is the basis dual to $\{S_c\}$. 
It is worth noting that the $\kappa(K,c)$ are the colored Kauffman brackets
of the knot (which  are  a version of the colored 
Jones polynomials of the knot). Indeed, the $c$-th coefficient of the
series expansion is  computed by plugging $S_c$ along the core
of the  regular neighborhood
of the knot and evaluating in $S^3$, that is by coloring the
knot with the $c$-th Jones-Wenzl idempotent and evaluating
the result in the skein space of the plane.

Let $\mathcal{F}_t(K)$ be the annihilator of $ \hat{Z}(K)$ in $K_t(T^2 \times I)$.
This is a left ideal that we call the {\em formal ideal} of $K$. 

\begin{theorem} The ideal $\mathcal{I}_t(K)$ lies in $\mathcal{F}_t(K)$.\end{theorem}

\proof Recall that a skein $\alpha$ in $T^2 \times I$ is in ${\mathcal I}_t(K)$ if when you glue
the $0$ end of $T^2 \times I$ to $X$, the skein is equivalent to $0$ in
$K_t(X)$. 
Let $\alpha \in {\mathcal I}_t(K)$ and $x\in K_t(S^1\times K^2)$. Then
$\alpha Z(K)$ is a functional on $K_t(S^1\times D^2)$ and its value
on $x$ is computed by embedding the skein $a$ in the knot
complement by gluing the cylinder over the torus to the 
knot complement, and then gluing the solid torus
with the skein $x$ in it to the knot complement. 
But the skein $\alpha$ can be transformed into the zero skein
by skein moves taking place entirely in the complement
of $x$, hence the value of the functional is zero.
Thus the functional itself is zero. 
\qed


If $\alpha$ is a skein in $K_t(T^2\times I)$, then $\alpha$ has a matrix
representation coming from the action on the
skein module of the solid torus, with basis $\{S_c\}$. The 
theorem shows that the rows of this matrix are orthogonal
to the vector $\hat{Z}(K)$.

\section{ The view beyond} The first question about this set up is whether
 $\mathcal{I}_t(K) = \mathcal{F}_t(K)$. Even in restricted cases, like when
the colored Jones polynomials of the knot are all trivial this is a provocative
question. This gives the first glimmer of how one might go about proving that the
Jones polynomial distinguishes the unknot.

The fact that the peripheral ideal annihilates
the vector whose coordinates are the colored Kauffman brackets
of the knot gives rise to a family of relations for
the generators which  depend on the parameter $t$. Do these depend
smoothly on $t$? If yes, they can be differentiated
and evaluated at $t=-1$ to yield derivatives of the A-polynomial of a knot. How are
those to be understood in terms of the $\G$-representation of the knot? The pairing used
in the orthogonality relation is clearly some sort of averaging process. Are there
integral formulas for the colored Jones polynomials, where the integral is taken over
the character variety of the knot?

 What kind
of information   about incompressible surfaces in a knot complement can be extracted 
from the noncommutative invariant. Can you see the placement of incompressible surfaces
in the knot complement from the colored Jones polynomials of the knot?


\begin{thebibliography}{00000}
\bibitem[AL]{AL} William Adams, Phillippe Loustaunau, {\em An introduction to Gr\"{o}bner
bases}, Graduate studies in Mathematics, ISSN 1065-7339; 3, AMS, Providence RI, 1994.
\bibitem[B]{B} D. Bullock,{\em  Rings of $Sl_2(\mathbb{C})$-characters and the 
Kauffman bracket skein module}, preprint (1996).
\bibitem[BH]{BH} G. W. Brumfiel, H. M. Hilden, {\em $Sl(2)$ Representations of 
Finitely Presented Groups},  Contemp. Math. {\bf 187} (1995).
\bibitem[BP]{BP}  D. Bullock, J.H. Przytycki, {\em Kauffman bracket skein module 
quantization of symmetric algebra and $so(3)$}, preprint
\bibitem[C]{C} P.M. Cohn, {\em Free Rings and their Relations}, Academic Press, 1971
\bibitem[CCGLS]{CCGLS} D. Cooper, M. Culler, H. Gillett, D.D. Long, P.B. Shalen, {\em Plane Curves 
associated to character varieties of 3-manifolds}, Inventiones Math. {\bf  118}, pp. 47-84 (1994).
\bibitem[CL]{CL} D. Cooper, D. Long, {\em Representation Theory and the A-polynomial of a knot},
Chaos, Solitons and Fractals, {\bf 9} (1998) no 4/5, 749--763.
\bibitem[Co]{Co} A. Connes, {\em Noncommutative Geometry}, Academic Press, London, 1994.
\bibitem[CS]{CS} M. Culler and P.B. Shalen,
{\em Varieties of group representations and splittings of 3-manifolds},
Ann. of Math. {\bf 117} (1983), 109-146.
\bibitem[FG]{FG} C. Frohman, R. Gelca, {\em Skein Modules and the Noncommutative Torus},
Preprint,xxx.lanl.gov/math.QA-9806107.
\bibitem[K]{Kas} Christian Kassel  {\em Quantum Groups}, Springer Verlag, 1995.
\bibitem[KM]{KM} M. Kapovich, J.J. Millson, {\em On Representation Varieties of Artin groups,
projective arrangements and fundamental groups of smooth complex algebraic varieties},
Preprint(1997)
\bibitem[Li]{Li} W. B. R. Lickorish, {\em An Introduction to Knot Theory}, 
  Springer, GTM {\bf 175}, 1997.
\bibitem[LM]{L-M} A. Lubotzky, A. Magid, {\em Varieties of representations of 
finitely generated groups},  Memoirs of the AMS {\bf 336} (1985).
\bibitem[PS]{PS} J.H. Przytycki, A.S. Sikora, {\em On Skein Algebras And 
$Sl_2(\mathbb{C})$-Character Varieties}, Preprint (1997).
\bibitem[S1]{S} Adam Sikora, { \em A geometric method in the theory of $SL_n$-representations 
       of groups}, Preprin, xxx.lanl.gov/ math.RT-9806016, (1998). 
\bibitem[R1]{R1} J. Roberts {\em Skeins and Mapping Class Groups}, Math.Proc. Camb. Phil. Soc. 
{\bf 115}(1994)53-77.
\bibitem[R2]{R} J. Roberts, {\em Kirby Calculus in Manifolds with Boundary}, Preprint,(1996).
\bibitem[Ri]{Ri} M. Rieffel, {\em Deformation Quantization of Heisenberg Manifolds},
Commun. Math. Phys. {\bf 122} (1989), 531-562

\end{thebibliography}
\end{document}